\documentclass[11pt,a4paper]{article}
\usepackage{mathptmx} % rm & math
\usepackage[scaled=0.99]{helvet} % ss
\usepackage{courier} % tt
\usepackage[T1]{fontenc}
\usepackage{graphicx}
\usepackage{amsmath}
\usepackage{amsfonts}
\usepackage{amssymb}
\usepackage{eufrak}

\newcommand{\dint}{\displaystyle\int}

\newtheorem{definition}{Definition}[section]
\newtheorem{theorem}{Theorem}[section]
\newtheorem{proposition}{Proposition}[section]
\newtheorem{lemma}{Lemma}[section]
\newtheorem{remark}{Remark}[section]
\newtheorem{corollary}{Corollary}[section]

\def\HH{\EuFrak H}

\def \bop {\noindent\textit{Proof}}
\def \eop {\hbox{}\nobreak\hfill
\vrule width 2mm height 2mm depth 0mm
\par \goodbreak \smallskip}
\def \R{\mathbb{R}}

\def \E{\mathbb{E}}

\def \bf{\textbf}
\def \it{\textit}

\def \ni {\noindent}

\def \eop {\hbox{}\nobreak\hfill
 $\square$
\par \goodbreak \smallskip}

\begin{document}

\date{}
\title{On the $\frac 1H$-variation of the divergence integral with respect to fractional Brownian motion with Hurst parameter $H<\frac12$ }
\author{El Hassan Essaky \\  Universit\'{e} Cadi Ayyad\\
 Facult\'{e} Poly-disciplinaire\\
Laboratoire de Mod\'{e}lisation et Combinatoire\\
D\'{e}partement de Math\'{e}matiques et d'Informatique\\ BP 4162, Safi, Maroc \\ Email: essaky@uca.ma\vspace*{0.2in}\\
 David Nualart \thanks{D. Nualart is supported by the NSF grant DMS 1208625}  \\
  The University of Kansas\\ Department of Mathematics\\
  Lawrence, Kansas 66045, USA\\ Email: nualart@ku.edu}
\maketitle\maketitle

\begin{abstract}
In this paper, we study  the $\frac 1H$-variation of stochastic divergence integrals $X_t=\int_{0}^{t} u_{s} \delta B_{s}$ with respect to a fractional Brownian motion $B$ with Hurst parameter $H< \frac12$. Under suitable assumptions on the process $u$, we prove that the $\frac 1H$-variation of $X$ exists in $L^1(\Omega)$ and  is equal to $e_H \int_0^T |u_s|^{\frac{1}{H}} ds$, where $e_H = \E\left[|B_1|^{\frac{1}{H}}\right]$. In the second part of the paper, we establish an integral representation for the fractional Bessel Process $\|B_t\|$, where $B_t$ is a $d$-dimensional fractional Brownian motion with Hurst parameter $H<\frac 12$. Using a multidimensional version of the result on the $\frac 1H$-variation of divergence integrals, we prove that if $2dH^2>1$, then the divergence integral in the integral representation of the fractional Bessel process has a $\frac 1H$-variation equals to a multiple of the Lebesgue measure. 
\end{abstract}

\vskip0.2cm

\ni {\small {\bf{Key words:}}
 Fractional Brownian motion, Malliavin calculus, Skorohod integral, Fractional Bessel processes.}

 \vskip0.2cm


%\renewcommand{\thefootnote}{\fnsymbol{footnote}}
%
%\vskip1cm
\bigskip
\section{Introduction}
The fractional Brownian  motion (fBm for short) $B=\{B_{t} , t\in [0,T]\}$ with Hurst  parameter $H\in (0,1)$  is a Gaussian self-similar process with stationary increments.
This process was introduced by Kolmogorov \cite{kol} and studied by Mandelbrot and Van Ness in \cite{MN}, where a stochastic integral representation in terms of a standard
Brownian motion was established. The parameter
$H$ is called Hurst index from the
statistical analysis, developed by the climatologist Hurst \cite{hurst}. The self-similarity and stationary increments properties make the fBm an appropriate model for many applications in diverse fields from biology to finance. From the properties of the fBm, it follows that for every $\alpha >0$
$$
\E\left(|B_t-B_s|^{\alpha}\right) = \E\left(|B_1|^{\alpha}\right)|t-s|^{\alpha H}.
$$
As a consequence of the Kolmogorov continuity theorem, we deduce that there exists a version of the fBm $B$ which is a continuous process and whose paths are $\gamma$-H\"{o}lder continuous for every $\gamma <H$.
Therefore, the fBm with Hurst parameter $H\neq \frac12$ is not a semimartingale and then the It\^{o} approach to the construction of stochastic integrals with respect to fBm is not valid. Two main approaches have been used in the literature to define stochastic integrals with respect to fBm with Hurst parameter $H$. Pathwise  Riemann-Stieltjes  stochastic integrals can be defined using Young's integral \cite{young} in the case $H>\frac 12$. When $H\in (\frac14, \frac12)$, the rough path analysis introduced by Lyons \cite{lyons} is a suitable method to construct pathwise stochastic integrals.

A second approach to develop a stochastic calculus with respect to the fBm is based on the techniques of Malliavin calculus. The divergence operator, which is the adjoint of the derivative operator, can be regarded as a stochastic integral, which coincides with the limit of Riemann sums constructed using the Wick product.
 This idea has been developed by
Decreusefond and \"{U}st\"{u}nel \cite{DU}, Carmona, Coutin and Montseny \cite{CC}, Al\`os, Mazet and Nualart \cite{AMN1, AMN2}, Al\`os and Nualart \cite{AN} and Hu \cite{hu}, among others.  The integral constructed by this method has zero mean. 
Different versions of the  It\^o formula have been proved by the divergence integral in these papers. 
In particular, if $H\in (\frac 14, 1)$ and $f\in C^2(\R)$ is a real-valued function satisfying some suitable  growth condition,  then the stochastic process $\{f'(B_t)\bf{1}_{[0,t]}, 0\le t \le T\}$ belongs to domain of the divergence operator and 
\begin{equation}\label{ito}
f(B_t) = f(0) + \dint_0^t f'(B_s)\delta B_s
 + H\dint_0^tf''(B_s) s^{2H-1} ds.
\end{equation}
For $H\in (0, \frac 14]$, this formula still holds, if the stochastic integral is interpreted as an extended divergence operator (see \cite{CN,LN}). A multidimensional version of the change of variable formula for the divergence integral has been recently proved by Hu, Jolis and Tindel in \cite{HMS}.

Using the self-similarity of fBm and  the Ergodic Theorem  one can prove that the fBm has a finite $\frac 1H$-variation on any interval $[0,t]$, equals to $e_H t$, where $e_H = \E\left[|B_1|^{\frac{1}{H}}\right]$ (see, for instance, Rogers \cite{rogers}). More precisely, we have, as $n$ tends to infinity  
\begin{equation}\label{res1}
 \sum_{i=0}^{n-1}|B_{t(i+1)/n}-B_{it/n}| ^{\frac{1}{H}}   \overset{L^{1}(\Omega)} {\longrightarrow}
t\, e_H. 
\end{equation}
This result has been generalized  by Guerra and Nualart \cite{GN}  to the case of divergence integrals with respect to the fBm  with Hurst parameter $H\in (\frac12, 1)$.

The purpose of this paper  is to study the $\frac 1H$-variation of   divergence processes
$X=\{ X_t, t\in [0,T]\}$, where  $X_t=\int_{0}^{t} u_{s} \delta B_{s}$, with respect to the fBm with Hurst parameter $H< \frac12$. Our main result, Theorem \ref{the2}, states that the $\frac 1H$-variation of $X$ exists in $L^1(\Omega)$ and  is equal to $e_H \int_0^T |u_s|^{\frac{1}{H}} ds$, under suitable assumptions on the integrand $u$. This is done by proving an estimate of the $L^{p}$-norm of the Skorohod integral $\int_{a}^{b} u_{s} \delta B_{s}$, where $0\leq a\leq b\leq T$.    Unlike the case $H>\frac 12$, here we need to impose H\"older continuity conditions on the process $u$ and its Malliavin derivative.
We also derive an extension of this result to divergence integrals with respect to a  $d$-dimensional fBm, where $d\ge 1$.  

 In the last part of the paper, we study the  fractional Bessel process $R= \{R_t, t\in [0,T]\}$, defined by $R_t:= \|B_t\|$, where $B$ is a $d$-dimensional fractional Brownian motion with Hurst parameter $H<\frac 12$.
 The following  integral representation of this process
\begin{equation}  \label{rep1}
R_t =  \displaystyle\sum_{i=1}^{d}\int_ 0^t\dfrac{B_s^{(i)}}{R_s}\delta B_s^{(i)} + H(d-1)\int_ 0^t \dfrac{s^{2H-1}}{R_s}ds,
\end{equation}
has been derived in \cite{GN} when $H>\frac 12$. Completing the analysis initiated in \cite{HN}, we
establish the representation  (\ref{rep1}) in the case $H<\frac 12$, using a suitable notion of the extended domain of the divergence operator. Applying the results obtained in the first part of the paper and assuming $2dH^2>1$, we prove  that the $\frac 1H$-variation of the divergence integral of the process  
$$\Theta_t:= \displaystyle\sum_{i=1}^{d}\int_ 0^t\dfrac{B_s^{(i)}}{R_t}\delta B_s^{(i)},$$ 
exists in $L^1(\Omega)$ and is equal to$ \dint_{\R^d}\left[\dint_  0^T \left | \left\langle\dfrac{B_s}{R_s}, \xi\right \rangle \right|^{\frac{1}{H}}ds\right]\nu(d\xi),
$ where $\nu$ is the normal distribution $N(0, I)$ on $\R^d$. We also discuss  some other properties of the process $\{\Theta_t,t\in [0,T]\} $.

The paper is organized as follows. Section 2 contains some preliminaries on Malliavin calculus.    In Section 3, we prove an $L^p$-estimate for divergence integral with respect to fBm. Section 4, is devoted to the study the $\frac 1H $-variation of the divergence integral with respect to fBm, for $H<\frac12$. Section 5, deals with   the $\frac 1H $-variation of the divergence integral with respect to $d$-dimensional fBm.  An application to fractional Bessel process has been given in Section 6.

\section{Preliminaries on Malliavin calculus}
Here we describe the elements from stochastic analysis that we will need in the paper. Let $B=\{B_{t} , t\in [0,T]\}$ be a   fractional Brownian  motion with Hurst parameter $H\in (0,1)$ defined in
a complete probability space $(\Omega, \mathcal{F},P)$,
 where $\mathcal{F}$ is generated by $B$. That is, $B$ is a 
centred Gaussian process with covariance function
\begin{equation*}
R_H(t,s):=\E(B_{t}B_{s}) = \dfrac12 (t^{2H}+s^{2H}-|t-s|^{2H}),
\end{equation*}
for $s,t \in [0,T]$.
We denote by $\HH$   the  Hilbert space associated to $B$,  defined as the closure of the linear space    generated by
the indicator functions $\{ \mathbf{1}_{[0,t]}, t\in [0,T]\} $, with respect to the inner product
\begin{equation*}
\langle \mathbf{1}_{[0,t]} , \mathbf{1}_{[0,s] } \rangle _{\HH}
=R_H(t,s), \hskip0.5cm s,t\in [0,T].
\end{equation*}
The mapping $\mathbf{1}_{[0,t]} \to B_{t}$ can be extended to a linear isometry between $\HH$ and the  Gaussian space  generated by $B$. We
denote by $B(\varphi)= \int_0^T \varphi_t dB_t $ the image of an element $\varphi \in \HH$
by this isometry.

We will first introduce some elements of the Malliavin calculus associated
with $B$. We refer to \cite{nualart} for a detailed account of these notions.
For a smooth and cylindrical  random variable $F=f\left( B(\varphi _{1}), \ldots ,
B(\varphi_{n})\right) $, with $\varphi_{i} \in \HH$ and $f\in
C_{b}^{\infty}(\R^{n})$ ($f$ and all its partial derivatives are bounded), the
derivative of $F$   is the $\HH$-valued random variable defined by
\begin{equation*}
D F =\sum_{j=1}^{n}\frac{\partial f}{\partial x_{j}}(B(\varphi_{1}),\dots,B(%
\varphi_{n}))\varphi_{j}.
\end{equation*}
For any integer $k\ge 1$ and any real number $p\ge 1$ we denote by $\mathbb{D%
}^{k,p}$ the Sobolev space defined as the the closure of the space of smooth and cylindrical 
random variables with respect to the norm
\begin{equation*}
\Vert F\Vert_{k,p}^{p}=\E(|F|^{p})+\sum_{j=1}^{k}  \E (\Vert
D^{j}F\Vert_{ \HH ^{\otimes j} }^{p}).
\end{equation*}
Similarly, for a given Hilbert space $V$ we can define Sobolev spaces of $V$%
-valued random variables $\mathbb{D}^{k,p}(W)$.

The divergence operator $\delta$ is introduced as the adjoint of the derivative operator. More precisely, an element $u\in L^{2}(\Omega;\HH)$ belongs to the domain of $\delta$, denoted by ${\rm Dom}\, \delta$,  if there exists
a constant $c_u$ depending on $u$ such that
\begin{equation*}
|\E(\langle D F,u\rangle_{\HH})|\leq c_u\Vert F\Vert_{2},
\end{equation*}
for any smooth random variable $F\in \mathcal{S}$. For any $u\in {\rm Dom}\, \delta$, $\delta(u)$ is the
element of $L^{2}(\Omega)$ given by the duality relationship
\begin{equation*}
\E(\delta (u)F)=\E(\langle D F,u\rangle_{\HH}),
\end{equation*}
for any $F\in \mathbb{D}^{1,2}$. We will make use of the notation $\delta
(u)=\int_{0}^{T}u_{s}\delta B_{s}$, and we call $\delta(u)$ the divergence integral of $u$ with respect to the fBm $B$. 
Note that $\E(\delta (
u ) )=0$. On the other hand, the space $\mathbb{D}^{1,2}(\HH)$ is included in the domain of $\delta $, and for $u\in \mathbb{D}^{1,2}(\HH)$,  the variance of $\delta(u)$ is given by
\begin{equation*}
\E(\delta (u)^{2})=\E(\Vert u\Vert_{\HH}^{2})+\E(\langle D u,(D
u)^{\ast}\rangle_{\HH\otimes\HH} ), 
\end{equation*}
where $(D u)^{\ast}$ is the
adjoint of $D u$ in the Hilbert space $\HH\otimes\HH$.
By Meyer's inequalities (see Nualart \cite{nualart}), for all $p>1$, the divergence operator 
is continuous from $ \mathbb{D}^{1,p}(\HH)$ into $ L^p(\Omega)$, that is,
\begin{equation}\label{meyer}
\E(|\delta (u)|^{p})\leq C_{p}\left( \E(\Vert u\Vert_{\HH%
}^{p})+\E(\Vert D u\Vert_{\HH\otimes\HH}^{p})\right).
\end{equation}
 We will make use of the property
\begin{equation}\label{p1}
\delta (Fu)= F\delta (u)+\langle D F,u\rangle_{\HH},
\end{equation}
which holds if  $F\in \mathbb{D}^{1,2}$,  $u\in {\rm Dom} \, \delta$ and the right-hand side is square integrable. We have also the commutativity relationship between $%
D $ and $\delta $
\begin{equation*}  
D \delta (u)= u + \int_{0}^{T} D u_{s}\delta B_{s},
\end{equation*}
which holds if $u\in \mathbb{D}^{1,2}(\HH)$ and the $\HH$-valued process $\{D  u_s, s\in
[0,T]\}$ belongs to the domain of $\delta $.

 The covariance of the fractional Brownian motion can be written as
 $$
R_H(t,s) = \int_0^{t\wedge s} K_H(t,u)K_H(s,u)du,
$$
where $K_H(t,s)$ is a square integrable kernel, defined for $0<s<t<T$. In what follows, we assume that $0<H <\frac12$. In this case, this kernel has the following expression
$$
K_H(t,s)= c_H\left[ \left(\frac{t}{s}\right)^{H-\frac12}(t-s)^{H-\frac12} -(H-\frac12)s^{H-\frac12}\int_s^t u^{H-\frac32}(u-s)^{H-\frac12}du\right],
$$
with  $c_H = \left(\frac{2H}{(1-2H)\beta(1-2H, H+\frac12)}\right)^{\frac12}$ and $\beta(x,y):= \dint_ 0^1 t^{x-1}(1-t)^{y-1}dt$ for $x, y>0$.  Notice also that 
$$
\frac{\partial K_H}{\partial t}(t,s) = c_H (H-\frac12)\left(\frac{t}{s}\right)^{H-\frac12}(t-s)^{H-\frac32}.
$$
From these expressions it follows that the  kernel $K_H$ satisfies  the following two estimates
\begin{equation}\label{est1A}
\left|\frac{\partial K_H}{\partial t}(t,s)\right| \leq c_H (t-s)^{H-\frac32},
\end{equation}
and 
\begin{equation}\label{est2}
|K_H(t,s)|\leq  d_H \left((t-s)^{H-\frac12} + s^{H-\frac 12} \right),
\end{equation}
for some constant $d_H$. 

Let $\mathcal{E}$ be the linear span of the indicator functions on $[0,T]$.
Consider the linear operator $K_H^*$ from ${\mathcal{E}}$ to $L^2([0, T])$ defined by
\begin{equation}\label{est0}
K_H^*(\varphi)(s) = K_H(T,s)\varphi(s)+ \int_s^T (\varphi(t)-\varphi(s))\dfrac{\partial K_H}{\partial t}(t,s)dt.
\end{equation}
 Notice that 
$$
K_H^*(\bf{1}_{[0, t]})(s) = K_H(t,s)\bf{1}_{[0, t]}(s).
$$
The operator $K_H^*$ can be expressed in terms of fractional derivatives as follows 
$$
(K_H^*\varphi)(s)= c_H \Gamma(H+\frac12)s^{\frac12 -H}(D_{T-}^{\frac12 -H} u^{H -\frac12}\varphi(u))(s).
$$
In this expression,  $D_{t-}^{\frac 12 -H}$ denotes the left-sided fractional derivative operator, given by
$$
D_{t-}^{\frac 12 -H }f(s):= \frac{1}{\Gamma(\frac 12+H )}\left(\dfrac{f(t)}{(t-s)^{\frac 12-H}}+\left(\frac 12 -H\right)\dint_s^t\dfrac{f(s)-f(y)}{(y-s)^{\frac 32-H}}dy\right),
$$
for almost all $s\in (0,t)$ and for a function $f $ in
the image of $L^p([0,t])$, $p\ge 1$,  by the left-sided fractional operator $I^{\frac 12-H}_{t-}$ (see \cite{SK} for more details).  
As a consequence  $C^{\gamma}([0,T])\subset \HH\subset L^2([0,T])$. It should be noted that the operator $K_H^*$ is an isometry between the Hilbert space $\HH$ and $L^2([0,T])$. That is, for every $\varphi, \psi\in\HH$,
\begin{equation}   \label{equ1}
\langle  \varphi, \psi \rangle_\HH= \langle K_H^*\varphi, K_H^* \psi \rangle_{L^2([0,T])}.
\end{equation}

 Consider the following seminorm on the space ${\mathcal{E}}$ 
\begin{equation}\label{iso}
\begin{array}{ll}
\| \varphi\|_ K^2 = \dint_ 0^T &\varphi^2(s)[(T-s)^{2H-1}+ s^{2H-1}]ds \\ & + \dint_0^T\left(\dint_s^T |\varphi(t)-\varphi(s)|(t-s)^{H-\frac32}dt\right)^2 ds.
\end{array}
\end{equation}
We denote by $\HH_K$ the completion of ${\mathcal{E}}$ with respect to this seminorm.  From the estimates (\ref{est1A}) and (\ref{est2}), there exists a constant $k_H$ such that for any $\varphi \in \HH_K$,
\begin{equation}\label{est01}
\| \varphi\|^2_{\HH} =\|K^*_{H}(\varphi)\|^2_{L^2([0,T])}
\leq k_H\| \varphi\|^2_{ K} .
\end{equation}
As a consequence, the space $\HH_K$  is continuously embedded in $\HH$.  This implies also that
$\mathbb{D}^{1, 2}(\HH_K) \subset \mathbb{D}^{1,2}(\HH) \subset {\rm Dom}\, \delta$.

One can show also that $\HH = I_{T-}^{\frac12 -H}(L^2([0,T]))$ (see \cite{DU}). Then, the space $\HH$ is too small for some purposes. For instance, it has been proved in \cite{CN}, that the trajectories of the fBm $B$
belongs to  $\HH$ if and only if $H>\frac14$.  This creates difficulties when defining the divergence $\delta(u)$ of a stochastic process whose trajectories do not belong to $\HH$, for example, if $u_t=f(B_t)$ and $H<\frac 14$, because the domain of $\delta$ is included in
$L^{2}(\Omega; \HH)$.  To overcome this difficulty, an extended domain of the divergence operator has been introduced in  \cite{CN}.  The main ingredient in the definition of this extended domain is the extension of the inner produce $\langle  \varphi, \psi \rangle_\HH$ to the case where  $\psi \in \mathcal{E}$ and $\varphi \in L^\beta([0,T])$ for some $\beta >\frac 1{2H}$ (see \cite{LN}).
More precisely, for  $\varphi \in L^\beta([0,T])$ and $\psi = \sum_{j=1}^{m}b_j\bf{1}_{[0,t_j]} \in \mathcal{E}$ we set
\begin{equation}  \label{ext}
\langle\varphi, \psi \rangle_\HH = \displaystyle\sum_{j=1}^{m}b_j\dint_0^T\varphi_s \dfrac{\partial R}{\partial s}(s, t_j)ds.
\end{equation}
 This expression coincides with the inner produce in $\HH$ if $\varphi \in \HH$, and it is well defined, because
\[
|\langle\varphi, \bf{1}_{[0,t]} \rangle_\HH|
 =   \left|\int_0^T\varphi_s \dfrac{\partial R}{\partial s}(s, t)ds \right|
   \leq \|\varphi\|_{L^\beta([0,T])} \sup_{0\leq t\leq T} \left(\int_0^T|\dfrac{\partial R}{\partial s}(s, t_j)|^{\alpha}ds\right)^{\frac{1}{\alpha}}<\infty.
\]

\vspace{0.1cm}
We will make use of following notations: for each $(a, b)\in\R^2$, $a\wedge b = \min(a, b)$ and
$a\vee b = \max(a, b)$.

\section{$L^p$-estimate of  divergence integrals with respect to fBm}

Let $V$ be a given Hilbert space. We introduce the following hypothesis for  a $V$-valued     stochastic process $u=\{ u_t, t\in [0,T]\}$, for some $p\ge 2$.

\medskip
\noindent
\textbf{Hypothesis}  $\mathbf{(A.1)}_p$  \textit{ Let $p\ge 2$. Then,   $\displaystyle\sup_{0\leq s\leq T}\Vert u_s\Vert_{L^{p}(\Omega; V)}  <\infty $ and there exist constants $L>0$, $0<\alpha <\frac12$ and $\gamma >\frac12 -H$ such  that,   
\begin{equation*}  \label{A1}
\Vert u_t -u_s\Vert_{L^{p}(\Omega; V)}\leq L s^{-\alpha }|t-s|^{\gamma},
\end{equation*}
  for all $0<s\leq t \leq T$. }

For any $ 0\le a< b \le T$, we will make use of the notation 
\[
\|u\|_{p,a,b} = \sup_{a\le s\le b} \|u_s\| _{L^{p}(\Omega; V)}.
\]
 The following lemma is a crucial ingredient to establish the 
  $L^p$-estimates for the divergence integral with respect to fBm.  
  
\begin{lemma}\label{lem1}
Let $u=\{u_t, 0\leq t\leq T\}$ be a process with values in a Hilbert space $V$, satisfying assumption  $\mathbf{(A.1)}_p$  for some $p\geq 2$. Then, there exists a positive constant $C$ depending on $H$, $\gamma$  and $p$ such that for every $0<a\leq b \le T$
\begin{equation}  \label{est1}
\E \left( \|  u  {\mathbf 1}_{[a,b]}   \| ^p_{\HH \otimes V} \right)   \leq C\left(\| u\|_{p,a,b}^p(b-a)^{pH}+ L^pa^{-p\alpha }(b-a)^{p\gamma +pH}\right).
\end{equation}
Moreover if $a=0$, then 
\begin{equation}  \label{est1a}
\E \left(  \|u {\mathbf 1}_{[0,b]}  \|^p_{\HH \otimes V}   \right) \leq  C\left(\| u\|_{p,0,b}^p b^{pH}+L^pb^{-p\alpha +p\gamma+pH}\right).
\end{equation}
\end{lemma}

\bop. Suppose first that $a> 0$. By equalities (\ref{equ1}) and (\ref{est0}) we obtain
\begin{eqnarray*}
&& \E \left( \|  u  {\mathbf 1}_{[a,b]}  \| ^p_{\HH \otimes V} \right) 
  = \E  \left( \| K_H^*(u \bf{1}_{[a,b]} ) \| ^p_{L^2([0,T];V)}  \right)  \\
 & &= \E \left(  \left\| K_H(T,s)u _s{\mathbf 1}_{[a,b]}(s)  +\dint_{s}^T\Big(u_t{\mathbf 1}_{[a,b]}(t)-u_{s} {\mathbf 1}_{[a,b]}(s)\Big)\dfrac{\partial K_H}{\partial t}(t,s)dt  \right \|^p_{L^2([0,T];V)}   \right).
\end{eqnarray*}
 Consider the decomposition 
\begin{eqnarray*}
 &&\dint_s^T\Big(u_t{\mathbf 1}_{[a,b]}(t)-u_s{\mathbf 1}_{[a,b]}(s)\Big)\dfrac{\partial K_H}{\partial t}(t,s)dt  
  = \left[\dint_s^b(u_t -u_s)\dfrac{\partial K_H}{\partial t}(t,s)dt\right]{\mathbf 1}_{[a,b]}(s)  \\
  &&\qquad +\left[-\dint_b^T u_s\dfrac{\partial K_H}{\partial t}(t,s)dt\right]{\mathbf 1}_{[a,b]}(s)    
 +\left[ \dint_a^b u_t\dfrac{\partial K_H}{\partial t}(t,s)dt\right]{\mathbf 1}_{[0,a]}(s)  \\
&& \qquad  := I_1 +I_2+I_3.
\end{eqnarray*}
Therefore
\[
   \E \left( \|  u  {\mathbf 1}_{[a,b]}  \| ^p_{\HH \otimes V} \right)    \le  C\sum_{i=0}^3 A_i,
   \]
   where $A_0= \E\left[\| K_H(T,\cdot)u \bf{1}_{[a,b]} \|^p_{L^2([0,T]; V)} \right]$  and for $i=1,2,3$,  $A_i= \E \left[\| I_i\|^p_{L^2([0,T];V)} \right]$.
Let us now estimate the four terms  $A_i$, $i=0,1,2,3$, in the previous inequality. By estimate (\ref{est2}),  Minkowski inequality and  Hypothesis $\mathbf{(A.1)}_p$ we obtain
\begin{eqnarray}\notag
A_0
 & \leq & C \E\left(\dint_a^b [(T-s)^{2H-1}+ s^{2H-1}]\Vert u_s\Vert^2_Vds\right)^{\frac{p}{2}} \\ \notag
 &\leq &   C\left(\dint_a^b [(T-s)^{2H-1}+s^{2H-1}]\| u_s\|^2_{{L^{p}(\Omega; V)}}ds\right)^{\frac{p}{2}} \\
 & \leq & C \| u \|_{p,a,b}^p (b-a)^{pH},\label{eqA0}
\end{eqnarray}
where we have used that $(T-a)^{2H}\leq (T-b)^{2H} +(b-a)^{2H}$ and $b^{2H} -a^{2H} \le (b-a)^{2H}$.
Using Minkowski inequality, Hypothesis $\mathbf{(A.1)}_p$ and  estimate (\ref{est1A}), it follows that 
\begin{eqnarray}
A_1  \notag
  & \leq & \left(\dint_a^b \left\Vert \dint_s^b(u_t -u_s)\dfrac{\partial K_H}{\partial t}(t,s)dt\right\Vert^2_{{L^{p}(\Omega; V)}}ds\right)^{\frac{p}{2}}   \\  \notag
  & \leq &  \left(\dint_a^b \left( \dint_s^b\Vert u_t -u_s\Vert_{{L^{p}(\Omega; V)}}\left|\dfrac{\partial K_H}{\partial t}(t,s)\right|dt\right)^2ds
\right)^{\frac{p}{2}}  \\
 &\leq  & C L^p\left(\dint_a^b \left( \dint_s^b s^{-\alpha }(t-s)^{\gamma+H-\frac32}dt\right)^2ds\right)^{\frac{p}{2}}.
 \label{eqA00}  
\end{eqnarray}
We have
\begin{eqnarray*}\notag
  \dint_a^b \left( \dint_s^b s^{-\alpha}(t-s)^{\gamma+H-\frac32}dt\right)^2ds \notag
 &= &\dfrac{1}{(\gamma +H-\frac12)^2}\dint_a^bs^{-2\alpha }(b-s)^{2\gamma +2H-1}ds \\ \notag
  &\le &
\dfrac{1}{(\gamma +H-\frac12)^2 (2\gamma +2H)} a^{-2\alpha } (b-a)^{2\gamma +2H}. 
\end{eqnarray*}
Substituting this expression into inequality (\ref{eqA00}),  yields
\begin{equation}\label{eqA1}
A_1
 \leq 
C L^p  a^{-p\alpha } (b-a)^{p\gamma +pH}.
\end{equation}
By the same arguments as above, it follows from Minkowski inequality and Hypothesis $\mathbf{(A.1)}_p$ that
\begin{eqnarray} 
A_2
  & =&   \notag
 \left(\dint_a^b \left\Vert\dint_b^T u_s\dfrac{\partial K_H}{\partial t}(t,s)dt\right\Vert^2_{{L^{p}(\Omega; V)}}ds\right)^{\frac{p}{2}}  \\ \notag
 & \leq  & C \left(\dint_a^b\left(\dint_b^T \Vert u_s\Vert_{{L^{p}(\Omega; V)}}(t-s)^{H-\frac32}dt\right)^2ds\right)^{\frac{p}{2}} \\ \notag
  &\leq  & C\|u \|_{p,a,b}^p \left(\dint_a^b\left((T-s)^{H-\frac12}-(b-s)^{H-\frac12}\right)^2ds\right)^{\frac{p}{2}}    \\ \notag
 &\leq  &C\|u \|_{p,a,b}^p\Big((T-a)^{2H}-(T-b)^{2H})+ (b-a)^{2H}\Big)^{\frac{p}{2}}   \\    \label{eqA2}
 &\leq  & C\| u \|_{p,a,b}^p(b-a)^{pH},
\end{eqnarray}
where we have used that $(T-a)^{2H}-(T-b)^{2H} \leq (b-a)^{2H}$.\\
Finally, for the term $A_3$, we obtain in the same way  
\begin{eqnarray}\notag
A_3  \notag
 &\leq & \left(\dint_0^a  \left(\dint_a^b \Vert u_t\Vert_{{L^{p}(\Omega; V)}}|\dfrac{\partial K_H}{\partial t}(t,s)|dt\right)^2ds\right)^{\frac{p}{2}}   \\
 & \leq  &  C\| u \|_{p,a,b}^p\left(\dint_0^a  \left(\dint_a^b (t-s)^{H-\frac{3}{2}}dt\right)^2ds\right)^{\frac{p}{2}}  \notag \\
& \le & C\| u \|_{p,a,b}^p\left(\dint_0^a  \left((a-s)^{H-\frac12} -(b-s)^{H-\frac12}\right)^2ds\right)^{\frac{p}{2}} \notag \\
& \le & C \| u \|_{p,a,b}^p (b-a) ^{pH}.  \label{eqA3}
\end{eqnarray}
For the last inequality we have used the following computations
\begin{eqnarray*}
& & \dint_0^a  \left((a-s)^{H-\frac12} -(b-s)^{H-\frac12}\right)^2ds  \\
&& \quad =    \frac 1{2H} \left( b^{2H} + a^{2H} -(b-a) ^{2H} \right)
 -2\dint_0^a  (a-s)^{H-\frac12}(b-s)^{H-\frac12}ds   \\
 & & \quad  \leq    \frac 1{2H} \left( b^{2H} + a^{2H} -(b-a) ^{2H} \right)  -2\dint_0^a  (b-s)^{2H-1}ds   \\
 &&  \quad  \leq    \frac 1{2H} \left( (b-a) ^{2H} - (b^{2H} -a^{2H}) \right) \le \frac 1{2H} (b-a)^{2H}.  
\end{eqnarray*}
 The inequality (\ref{est1})  follows from the estimates (\ref{eqA0}), (\ref{eqA1}), (\ref{eqA2}) and (\ref{eqA3}). The case $a=0$ can be proved using  similar arguments. The proof of Lemma \ref{lem1} is then completed.     
 \eop
 \vspace{0.4cm}
We are now in the position to prove the following theorem which gives an estimate of the $L^{p}$-norm of the Skorohod integral of a process $u$ with respect to a fBm with Hurst parameter $H\in (0,\frac 12)$.  We first need the following assumption on the process $u$.

\medskip
\noindent
\textbf{Hypothesis}  $\mathbf{(A.2)}_p$  \textit{ Let  $u\in \mathbb{D}^{1, 2}(\HH)$ be a real-valued stochastic process, which satisfies Hypothesis  $\mathbf{(A.1)}_p$   with constants  $L_u$, $ \alpha_1$ and $\gamma$ for a fixed $p\geq 2$. We also assume that the $\HH$-valued process $\{Du_s, s\in [0,T]\}$ satisfies Hypothesis $\mathbf{(A.1)}_p$  with constants  $L_{Du}$, $ \alpha_2$ and $\gamma$ for the same  value of $p$.  }

\medskip  Hypothesis $\mathbf{(A.2)}_p$   means that  $u_s$ and $Du_s$ have bounded $L^p$ norms in $[0,T]$ and   satisfy
\begin{eqnarray}   
  \Vert u_t -u_s\Vert_{L^{p}(\Omega)}&\leq& L_us^{-\alpha _1}|t-s|^{\gamma} \label{assump1}
  \\
  \Vert  Du_t -Du_s\Vert_{L^{p}(\Omega; \HH)}&\leq & L_{Du}s^{-\alpha _2}|t-s|^{\gamma}, \label{assump2}
\end{eqnarray}
for all $0<s\leq t\leq T$.

\medskip

\begin{theorem}\label{the1}
Suppose that $u\in \mathbb{D}^{1, 2}(\HH)$ is a stochastic process  satisfying Hypothesis $\mathbf{(A.2)}_p$   for some $p\geq 2$. Let $0< a\leq b\leq T$.   Then, there exists a positive constant $C$  depending on $H$, $\gamma$ and $p$ such that 
 \begin{eqnarray}  \notag
& & \E\left( \left |\dint_a^b u_s \delta B_s\right|^p \right) \\ 
& \leq  &
C\left((\| u \|_{p,a,b}^p+\| Du \|_{p,a,b}^p)(b-a)^{pH}+ (L_u ^pa^{-p\alpha_1}+L_{Du}^pa^{-p\alpha_2})(b-a)^{p\gamma +pH}\right).\qquad
\label{ineq1}
\end{eqnarray}
If $a=0$, then
 \begin{eqnarray} 
 \E  \left( \left|\dint_0^b u_s \delta B_s\right|^p \right) \leq C\left((\|u \|_{p,a,b}^p+\|Du \|_{p,a,b}^p) b^{pH}+(L_u^pb^{-p\alpha_1}+L_{Du}^pb^{-p\alpha_2})b^{p\gamma+pH}\right).
\label{ineq2}
\end{eqnarray}
\end{theorem}
\bop. 
By inequality (\ref{meyer}), we have 
$$
\E \left( \left|\dint_a^b u_s \delta B_s\right|^p \right) \leq C_p\left ( \E (\|  u  {\mathbf 1}_{[a,b]}  \| ^p_{  \HH } )+\E ( \|   D_s(u_t\bf{1}_{[a,b]}(t)) \| ^p_{ \HH \otimes \HH}\right).
$$ 
The first and the second terms of the above inequality can be estimated applying   Lemma \ref{lem1}  to the processes $u$ and $Du$, with $V= \mathbb{R}$ and $V=\HH$, respectively.       Theorem \ref{the1} is then proved.\eop

\begin{remark} If we suppose  that $\alpha_1= \alpha_2=0$ in Hypothesis $\mathbf{(A.2)}_p$,  that is, $u$ and $Du$ are H\"older continuous in $L^p$ on $[0,T]$, then estimate  (\ref{ineq1}) in Theorem \ref{the1} can be written as
$$
\E\left( \left|  \dint_a^b u_s \delta B_s\right|^p \right) \leq C\Vert u\Vert_{1,p,\gamma}^p(b-a)^{pH},
$$
where
$$
\Vert u\Vert_{1,p,\gamma}=\displaystyle\sup_{0\leq s<t\leq T}\dfrac{\Vert u_t -u_s\Vert_{1,p}}{|t-s|^{\gamma}}+\displaystyle\sup_{0\leq s\leq T} \Vert u_s\Vert_{1,p}.
$$
\end{remark}

\section{The $\frac{1}{H}$-variation of divergence integral with respect to fBm}
Fix $q\geq 1$ and $T>0$ and set $t_i^n:= \frac{iT}{n}$, where $n$ is a positive integer and $i=0,1,2,\dots,n$. We need the following definition.
\begin{definition}
Let $X$ be a given stochastic process defined in the complete probability space $(\Omega, {\cal F}, P)$. Let $V_n^q(X)$ be the random variable defined by
$$
V_n^q(X):= \sum_{i=0}^{n-1}|\Delta_i^n X|^q,
$$  
where $\Delta_i^n X := X_{t^n_{i+1}}-X_{t^n_{i}}$. We define the $q$-variation of $X$ as the limit in $L^1(\Omega)$, as  $n$ goes to infinity, of $V_n^q(X)$ if this limit exists.
\end{definition}
 
As in the last section we assume that  $H\in (0, \frac12)$. In this section, we need the following assumption on the process $u$.

\medskip
\noindent
\textbf{Hypothesis}  $\mathbf{(A.3)}$   \textit{
  Let  $u\in \mathbb{D}^{1, 2}(\HH)$ be a real-valued stochastic process which is bounded in $L^q(\Omega)$ for some $q>\frac{1}{H}$ and satisfies  the  H\"older continuity property (\ref{assump1}) with $p=\frac{1}{H}$, that is 
\begin{eqnarray}   
  \Vert u_t -u_s\Vert_{L^{\frac{1}{H}}(\Omega)}&\leq& L_us^{-\alpha _1}|t-s|^{\gamma}. \label{assump11}
\end{eqnarray}
Suppose also that the $\HH$-valued process $\{Du_s, s\in [0,T]\}$  is bounded in $L^{\frac{1}{H}}(\Omega; \HH)$ and satisfies  the  H\"older continuity property  (\ref{assump2}) with $p=\frac{1}{H}$, that is 
\begin{eqnarray}   
  \Vert  Du_t -Du_s\Vert_{L^{\frac{1}{H}}(\Omega; \HH)}&\leq& L_{Du}s^{-\alpha _2}|t-s|^{\gamma}. \label{assump21}
\end{eqnarray}
Moreover, we assume that the derivative $\{D_tu_s, s,t\in [0,T]\}$ satisfies
\begin{equation}\label{assump3}
\displaystyle\sup_{0\leq s\leq T}\Vert D_su_t\Vert_{L^{\frac{1}{H}}(\Omega)} \leq K t^{-\alpha_3},
\end{equation}
for every $t\in(0, T]$ and for some constants $0<\alpha_3<2H$ and $K>0$.}

\medskip
Consider the indefinite divergence integral of $u$ with respect to the   fBm $B$, given by 
\begin{equation}  \label{equ2}
X_t = \int_0^t u_s \delta B_s := \delta(u\bf{1}_{[0,t]}).
\end{equation}
The main result of this section is the following theorem.

\begin{theorem}\label{the2}
Suppose that $u\in \mathbb{D}^{1, 2}(\HH)$ is a stochastic process  satisfying Hypothesis $\bf{(A.3)}$, and consider the divergence integral process $X$ given by (\ref{equ2}).  Then, we have 
\[
V_n^{\frac{1}{H}}(X) \overset{L^{1}(\Omega)} {\longrightarrow} e_H \dint_  0^T|u_s|^{\frac{1}{H}}ds, 
\]
as $n$ tends to infinity,
where $e_H = \E \left[|B_1|^{\frac{1}{H}}\right]$.
\end{theorem}
\bop. 
We need to show that the expression
\[
F_n:= \E\left(\left|\sum_{i=0}^{n-1}\left|\dint_{t_i^n}^{t_{i+1}^n} u_s \delta B_s\right|^{\frac{1}{H}}-e_H \dint_0^T |u_s|^{\frac{1}{H}}ds\right|\right),
\]
converges to zero as $n$ tends to infinity.
Using (\ref{p1}), we can write
\begin{equation}\label{decom}
\begin{array}{ll}
\dint_{t_i^n}^{t_{i+1}^n} u_s \delta B_s
 &=\dint_{t_i^n}^{t_{i+1}^n} (u_s-u_{t_i^n}) \delta B_s + \dint_{t_i^n}^{t_{i+1}^n} u_{t_{i}^n}\delta B_s
\\ & =\dint_{t_i^n}^{t_{i+1}^n} (u_s-u_{t_i^n}) \delta B_s-\langle Du_{t_i^n}, \bf{1}_{[t_{i}^n, t_{i+1}^n]}\rangle_{{\HH}} + u_{t_{i}^n}(B_{t_{i+1}^n}-B_{t_{i}^n}).
\\ & := A_{i}^{1,n} -A_{i}^{2,n} +A_{i}^{3,n}.
\end{array}
\end{equation}
By the triangular inequality, we obtain
\begin{equation}
 F_n \le   \E\left(\sum_{i=0}^{n-1}\left|  |A_{i}^{1,n} -A_{i}^{2,n} +A_{i}^{3,n} |^{\frac{1}{H}} - |A_{i}^{3,n} |^{\frac{1}{H}}\right|\right)  +D_n, \label{eq45}
\end{equation}
where
\[
D_n=\E\left(\left|
\sum_{i=0}^{n-1} |A_{i}^{3,n} |^{\frac{1}{H}}-e_H \dint_0^T |u_s|^{\frac{1}{H}}ds\right|\right).
\]
Using the mean value theorem and H\"older inequality, we can write
\begin{eqnarray}   \notag
& &\E\left(  \sum_{i=0}^{n-1}\left|  |A_{i}^{1,n} -A_{i}^{2,n} +A_{i}^{3,n}  |^{\frac{1}{H}} - |A_{i}^{3,n} |^{\frac{1}{H}}\right|\right)   \\  \notag
&& \leq  \frac{1}{H}\E\left(\sum_{i=0}^{n-1}  |A_{i}^{1,n} -A_{i}^{2,n} |\left[ |A_{i}^{1,n} -A_{i}^{2,n}+A_{i}^{3,n}|^{\frac{1}{H}-1} + |A_{i}^{3,n}|^{\frac{1}{H}-1}\right]\right)   \\   \notag
& & \leq 
C\left[  \E\left(\sum_{i=0}^{n-1} |A_{i}^{1,n} -A_{i}^{2,n} |^{\frac{1}{H}}\right)\right]^H \\
&& \qquad \qquad \times \left[\E\left(\sum_{i=0}^{n-1} |A_{i}^{1,n} -A_{i}^{2,n}+A_{i}^{3,n} |^{\frac{1}{H}}\right)
+\E\left(\sum_{i=0}^{n-1} |A_{i}^{3,n}  |^{\frac{1}{H}}\right)\right]^{1-H}.  \label{eq451}
\end{eqnarray}
Substituting (\ref{eq451}) into  (\ref{eq45}) yields
\[
F_n  \le CA_{n}^H(B_n + C_n)^{1-H} + D_n,
\]
where
\begin{eqnarray*}
A_n&=&\E\left(\sum_{i=0}^{n-1}  |A_{i}^{1,n} -A_{i}^{2,n} |^{\frac{1}{H}}\right), \\
B_n &=& \E\left(\sum_{i=0}^{n-1} |A_{i}^{1,n} -A_{i}^{2,n}+A_{i}^{3,n} |^{\frac{1}{H}}\right),\\
C_n &=&\E\left(\sum_{i=0}^{n-1} |A_{i}^{3,n}  |^{\frac{1}{H}}\right).
\end{eqnarray*}
The proof will be divided into several steps. Along the proof, $C$ will denote a generic constant, which may vary from line to line and  may depend on the processes $u$ and $Du$ and the different parameters appearing in the computations, but it is independent of $n$.
\\
 \it{Step 1.} We first prove that   $B_n$ and $C_n$ are bounded. Remark that
\begin{eqnarray*}
B_n &= &\E\left( \left|  \int_{0}^{\frac{T}{n}} u_s \delta B_s\right|^{\frac{1}{H}} \right)+ \E \left( \sum_{i=1}^{n-1}\left|\int_{t_i^n}^{t_{i+1}^n} u_s\delta B_s\right|^{\frac{1}{H}}  \right) \\
& := &K_1^n+K_2^n.
\end{eqnarray*} 
Using estimate (\ref{ineq2}) with $p=\frac{1}{H}$, it follows that 
\begin{eqnarray*} 
 K_1^n 
 &\leq & C\left(\| u\|  ^\frac{1}{H}_{\frac{1}{H},0,\frac{T}{n}}+\| Du\|^\frac{1}{H}_{\frac{1}{H},0,\frac{T}{n}}\right)n^{-1}+ \left(L_{u}^{\frac{1}{H}} n^{\frac{\alpha_1}{H}}+L_{Du}^{\frac{1}{H}}n^{\frac{\alpha_2}{H}}\right)n^{-\frac{\gamma}{H}-1}
\\ & \leq &  C \left(n^{-1} +n^{\frac{\alpha_1}{H}-\frac{\gamma}{H}-1}+n^{\frac{\alpha_2}{H}-\frac{\gamma}{H}-1}\right).
\end{eqnarray*}
 Therefore, $K_1^n$ is bounded since $\alpha_1 <\gamma +H$ and $\alpha_2 <\gamma +H$. In a similar way,  estimate (\ref{ineq1}) leads to
\begin{eqnarray*}
 K_2^n &\leq & C\sum_{i=1}^{n-1}\bigg\{\left(\| u\|^{\frac{1}{H}}_{\frac{1}{H},t_i^n,t_{i+1}^n}+\| Du \|^{\frac{1}{H}}_{\frac{1}{H},t_i^n,t_{i+1}^n}\right)(t_{i+1}^n-t_{i}^n) \\ 
 && \qquad\qquad\quad+ \left(L_{u}^{\frac{1}{H}}(t_i^n)^{-\frac{\alpha_1}{H}}+L_{Du}^{\frac{1}{H}}(t_{i}^n)^{-\frac{\alpha_2}{H}}\right)(t_{i+1}^n-t_{i}^n)^{\frac{\gamma}{H} +1}\bigg\}    \\
 & \leq  &C \left(1+ n^{\frac{\alpha_1}{H}-\frac{\gamma}{H}-1}\displaystyle\sum_{i=1}^{n-1}{i^{-\frac{\alpha_1}{H}}}+ n^{\frac{\alpha_2}{H}-\frac{\gamma}{H}-1}\displaystyle\sum_{i=1}^{n-1}{i^{-\frac{\alpha_2}{H}}}\right).
\end{eqnarray*} 
This proves that $K_2^n$ is bounded  and so is $B_n$.  
 Using H\"{o}lder inequality and the fact that $u$ is bounded in $L^q(\Omega)$ for $q>\frac{1}{H}$, we obtain
\begin{eqnarray*}
C_n
&=&\sum_{i=0}^{n-1}\E\left( |u_{t_{i}^n} |^{\frac{1}{H}} |B_{t_{i+1}^n}-B_{t_{i}^n})|^{\frac{1}{H}}\right)
\\ & \leq  & \sum_{i=0}^{n-1}\left[ \E\left(|u_{t_{i}^n}|^{q}\right)\right] ^{\frac{1}{qH}}\left[ \E\left(|B_{t_{i+1}^n}-B_{t_{i}^n})|^{\frac{q}{qH -1}}\right)\right]^{1-\frac{1}{qH}}
\\ & 
\leq  & C   \sum_{i=0}^{n-1}(t_{i+1}^n-t_{i}^n) =CT,
\end{eqnarray*}
and this proves the boundedness of $C_n$.\\
\it{Step 2.} We prove that $A_n$ converges to zero.  Consider the decomposition
\[
\sum_{i=0}^{n-1}|A_{i}^{1,n}| ^{\frac{1}{H}} 
 =\left|\int_{0}^{\frac{T}{n}} (u_s-u_{0}) \delta B_s\right|^{\frac{1}{H}} + \sum_{i=1}^{n-1}|A_{i}^{1,n} |^{\frac{1}{H}}.
\]
Using estimate (\ref{ineq2}) with $p =\frac{1}{H}$, it follows that 
\begin{eqnarray*}
&&  \E\left( \left|\int_{0}^{\frac{T}{n}} (u_s-u_{0}) \delta B_s \right|^{\frac{1}{H}}   \right) \\
&& \leq C\left[ \|u-u_{0}  \|^\frac{1}{H}_{\frac{1}{H},0,\frac{T}{n}}+ \| Du-Du_{0}  \|^\frac{1}{H}_{\frac{1}{H},0,\frac{T}{n}}\right]n^{-1}+ \left[L_{u}^{\frac{1}{H}}{n}^{\frac{\alpha_1}{H}}+L_{Du}^{\frac{1}{H}}{n}^{\frac{\alpha_2}{H}}\right]{n}^{-\frac{\gamma}{H}-1}  \\
&& \leq C n^{-1}\left(1  +{n^{\frac{\alpha_1}{H}-\frac{\gamma}{H}}}+{n^{\frac{\alpha_2}{H}-\frac{\gamma}{H}}}\right).
\end{eqnarray*}
Therefore $ \E\left( \left|\int_{0}^{\frac{T}{n}} (u_s-u_{0}) \delta B_s \right|^{\frac{1}{H}}   \right) $ converges to zero as $n$ tends to infinity, since $\alpha_1 <\gamma +H$ and $\alpha_2 <\gamma +H$. We can also prove that $\E \left(\sum_{i=1}^{n-1} |A_{i}^{1,n}|^{\frac{1}{H}}  \right)$ converges to zero.  In fact, using estimate (\ref{ineq1}) with $p =\frac{1}{H}$, we obtain 
\begin{eqnarray*}
\E\left( \sum_{i=1}^{n-1} |A_{i}^{1,n}|^{\frac{1}{H}}   \right)
 & \leq  &
 C \sum_{i=1}^{n-1}\Bigg[\left(\| u-u_{t_i^n} \|^{\frac{1}{H}}_{\frac{1}{H},t_i^n,t_{i+1}^n}+\| Du-Du_{t_i^n} \|^{\frac{1}{H}}_{\frac{1}{H},t_i^n,t_{i+1}^n}\right)(t_{i+1}^n-t_{i}^n)\\ 
 && \qquad\qquad\quad+ \left(L_{u}^{\frac{1}{H}}(t_i^n)^{-\frac{\alpha_1}{H}}+L_{Du}^{\frac{1}{H}}(t_{i}^n)^{-\frac{\alpha_2}{H}}\right)(t_{i+1}^n-t_{i}^n)^{\frac{\gamma}{H} +1}\Bigg]  \\
& & \leq C n^{-\frac{\gamma}{H}-1}\left({n^{\frac{\alpha_1}{H}}}\displaystyle\sum_{i=1}^{n-1}{i^{-\frac{\alpha_1}{H}}}+ {n^{\frac{\alpha_2}{H}}}\displaystyle\sum_{i=1}^{n-1}{i^{-\frac{\alpha_2}{H}}}\right), 
\end{eqnarray*}
where we have used the fact that
\begin{equation*}
\| u-u_{t_i^n} \|_{\frac{1}{H},t_i^n,t_{i+1}^n}   \leq L_{u}T^{\gamma-\alpha_1} i^{-\alpha_1} n^{\alpha_1-\gamma},
\end{equation*}
and 
\begin{equation*}
\| Du-Du_{t_i^n} \|_{\frac{1}{H},t_i^n,t_{i+1}^n} 
\leq L_{Du}T^{\gamma-\alpha_1} i^{-\alpha_2} n^{\alpha_2-\gamma}.
\end{equation*} 
From the above computations,\ it follows that $\E\left( \sum_{i=1}^{n-1}|A_{i}^{1,n}|^{\frac{1}{H}}   \right)$ converges to zero as $n$ goes to infinity. Therefore, we conclude that
\begin{equation}\label{term1}
\lim _{n \rightarrow \infty} \E\left( \sum_{i=1}^{n-1} |A_{i}^{1,n}|^{\frac{1}{H}}   \right) =0.
\end{equation}
Second, let us prove that $\E \left(\sum_{i=1}^{n-1}|A_{i}^{2,n}|^{\frac{1}{H}} \right)$ converge to zero as $n$ tends to infinity.  It follows from (\ref{ext}) that each term  $A_i^{2,n}$ can be expressed as 
\begin{equation*}
A_{i}^{2,n}
 =\int_0^T D_su_{t_i^n}\dfrac{\partial}{\partial s}\bigg(R(s,t_{i+1}^n)-(R(s,t_{i}^n)\bigg)ds.
\end{equation*}
Therefore we have the following decomposition
\begin{equation*}
A_{i}^{2,n}:= J_1^{i,n}+J_2^{i,n}+J_3^{i,n},
\end{equation*}
where
\begin{eqnarray*}
 J_1^{i,n} &=& \frac 12 
  \int_0^{t_{i}^n} D_su_{t_i^n}\dfrac{\partial}{\partial s}\left(((t_{i}^n-s)^{2H}-(t_{i+1}^n-s)^{2H}))\right)ds,  \\
  J_2^{i,n} &=&\frac 12 \int^{t_{i+1}^n}_{t_{i}^n} D_su_{t_i^n}\dfrac{\partial}{\partial s}\left(((s-t_{i}^n)^{2H}-(t_{i+1}^n-s)^{2H}))\right)ds, \\
   J_3^{i,n}&=& \frac 12 \int_{t_{i+1}^n}^{T} D_su_{t_i^n}\frac{\partial}{\partial s}\left(((s-t_{i}^n)^{2H}-(s-t_{i+1}^n)^{2H}))\right)ds.
\end{eqnarray*}
Using Minkowski inequality and assumption (\ref{assump3}), we obtain
\begin{eqnarray*}
\E \left(  \sum_{i=0}^{n-1}|J_1^{i,n}|^{\frac{1}{H}}   \right) 
 & \leq &  H \sum_{i=0}^{n-1}\left[ \int_0^{t_{i}^n} \Vert D_su_{t_i^n}\Vert_{L^{\frac{1}{H}}(\Omega)} \left|(t_{i+1}^n-s)^{2H-1}-(t_{i}^n-s)^{2H-1})\right|ds\right]^{\frac{1}{H}}   \\
& \leq  &  C \sum_{i=1}^{n-1}(t_{i}^n)^{-\frac{\alpha_3}{H}}\left[ \int_0^{t_{i}^n} \left[(t_{i}^n-s)^{2H-1}-(t_{i+1}^n-s)^{2H-1}\right]ds\right]^{\frac{1}{H}}  \\
& =& C \sum_{i=1}^{n-1}(t_{i}^n)^{-\frac{\alpha_3}{H}}\left[ (t_{i+1}^n-t_{i}^n)^{2H}-\left[(t_{i+1}^n)^{2H}-(t_{i}^n)^{2H}\right]\right]^{\frac{1}{H}}  \\
 & \leq  & C \sum_{i=1}^{n-1}(t_{i}^n)^{-\frac{\alpha_3}{H}}(t_{i+1}^n-t_{i}^n)^{2}  \\
 & \leq  & C {n^{\frac{\alpha_3}{H} -2}}\sum_{i=1}^{n-1} {i^{-\frac{\alpha_3}{H}}}.
\end{eqnarray*}
Taking into account that $\alpha_3 <2H$, we obtain that  $\E \left(  \sum_{i=0}^{n-1}|J_1^{i,n}|^{\frac{1}{H}}   \right) $ converges to zero as $n$ tends to infinity. By means of similar arguments, we can show that $\E \left(  \sum_{i=0}^{n-1}|J_2^{i,n}|^{\frac{1}{H}}   \right) $ and $\E \left(  \sum_{i=0}^{n-1}|J_3^{i,n}|^{\frac{1}{H}}   \right) $ converge to zero as $n$ tends to infinity. Therefore,
\begin{equation}\label{term2}
\lim _{n \rightarrow \infty} \E\left( \sum_{i=1}^{n-1}|A_{i}^{2,n}|^{\frac{1}{H}}   \right) =0.
\end{equation}
Consequently,   from (\ref{term1}) and (\ref{term2}) we deduce that  that $A_n$ converge to zero as $n$ goes to infinity.\\

\noindent
\it{Step 3.}  In order to show that the term $D_n$ converges to zero as $n$ tends to infinity, we replace $n$ by the product $nm$ and we let first $m$  tend to infinity.  That is, we consider the partition of interval $[0,T]$ given by $0=t_0^{nm}<\cdots <t_{nm}^{nm} = T$ and we define
\begin{eqnarray} 
Z^{n,m} &:=&   \notag
\left|  \sum_{i=0}^{nm-1} |u_{t_{i}^{nm}}|^{\frac{1}{H}}|\Delta_{i}^{nm}B|^{\frac{1}{H}}
 -e_H  \sum_{j=0}^{n-1} |u_{t_{j}^{n}}|^{\frac{1}{H}}(t_{j+1}^{n} -t_{j}^{n})\right|   \\  \notag
  & = & \bigg| \sum_{j=0}^{n-1} \bigg[\sum_{i=jm}^{(j+1)m -1} \left(|u_{t_{i}^{nm}}|^{\frac{1}{H}}-|u_{t_{j}^{n}}|^{\frac{1}{H}}\right)|\Delta_{i}^{nm}B|^{\frac{1}{H}}  \\   \notag
  & & \qquad +|u_{t_{j}^{n}}|^{\frac{1}{H}}\left(\sum_{i=jm}^{(j+1)m -1}|\Delta_{i}^{nm}B|^{\frac{1}{H}}-e_H (t_{j+1}^{n} -t_{j}^{n})\right)\bigg]\bigg|.   \\  \notag
 & \leq  &  \sum_{j=0}^{n-1} \sum_{i=jm}^{(j+1)m -1} \left||u_{t_{i}^{nm}}|^{\frac{1}{H}}-|u_{t_{j}^{n}}|^{\frac{1}{H}}\right||\Delta_{i}^{nm}B|^{\frac{1}{H}}  \\  \notag
   & &\qquad  +\sum_{j=0}^{n-1}|u_{t_{j}^{n}}|^{\frac{1}{H}}\left|\sum_{i=jm}^{(j+1)m -1}|\Delta_{i}^{nm}B|^{\frac{1}{H}}-e_H (t_{j+1}^{n} -t_{j}^{n})\right|    \\ \notag \label{eq3}
  & := & Z_1^{n,m} + Z_2^{n,m}.
\end{eqnarray}
By the mean value theorem, we can write
\[
Z^{n,m}_1  \le
 \frac 1H  \sum_{j=0}^{n-1} \sum_{i=jm}^{(j+1)m -1} \left|u_{t_{i}^{nm}}-u_{t_{j}^{n}}\right|\left(|u_{t_{i}^{nm}}|^{\frac{1}{H}-1}+|u_{t_{j}^{n}}|^{\frac{1}{H}-1}\right)|\Delta_{i}^{nm}B|^{\frac{1}{H}}.
\]
Using H\"{o}lder inequality, assumption (\ref{assump11}) as well as the boundedness of $u$ in $L^q(\Omega)$ for some $q>\frac{1}{H}$,  we obtain
\[
\E (Z^{n,m}_1)
\leq Cn^{-1}m^{-1} \sum_{j=0}^{n-1}  \sum_{i=jm}^{(j+1)m -1} (t_j^n)^{-{\alpha_1}}(t_i^{nm}-t_j^{n})^{{\gamma}}
\leq  {C}{n^{{\alpha_1}-{\gamma} -1}} \sum_{j=0}^{n-1} {j^{{-\alpha_ 1}}},
\]
which implies
\begin{equation}
\lim_{n\rightarrow \infty}\sup_{m\ge 1}\E (Z^{n,m}_1 )= 0.  \label{equ6}
\end{equation}
On the other hand, using H\"{o}lder inequality and the fact that $u$ is bounded in $L^q(\Omega)$ for some $q>\frac{1}{H}$, we have
\begin{eqnarray*}
\E(Z^{n,m}_2)
&\leq &
\sum_{j=0}^{n-1}\left[ \E(|u_{t_{j}^{n}}|^q)\right]^{\frac{1}{qH}}\left[ \E \left(\left|\sum_{i=jm}^{(j+1)m -1}|\Delta_{i}^{nm}B|^{\frac{1}{H}}-e_H (t_{j+1}^{n} -t_{j}^{n})\right|^{\frac{qH}{qH-1}}   \right)\right]^{1-\frac{1}{qH}}  \\
  &  \leq & C \sum_{j=0}^{n-1}\left[ \E \left( \left|\sum_{i=jm}^{(j+1)m -1}|\Delta_{i}^{nm}B|^{\frac{1}{H}}-e_H (t_{j+1}^{n} -t_{j}^{n})\right|^{\frac{qH}{qH-1}} \right)\right]^{1-\frac{1}{qH}}.
\end{eqnarray*}
For any fixed $n\ge 1$, by the Ergodic Theorem  the sequence $ \sum_{i=jm}^{(j+1)m -1}|\Delta_{i}^{nm}B|^{\frac{1}{H}}-e_H (t_{j+1}^{n} -t_{j}^{n})$ converges to $0$ in $L^s$ as $m$ tends to infinity, for every $s>1$. This implies that, for any $n\ge 1$,
\begin{equation}
\lim_{m\rightarrow \infty} \E (Z^{n,m}_2 )= 0.  \label{equ7}
\end{equation}
Therefore, it follows from (\ref{equ6}) and (\ref{equ7}) that 
\begin{equation}
\lim_{n\rightarrow \infty} \lim_{m\rightarrow \infty} \E (Z^{n,m})= 0.  \label{equ8}
\end{equation}
By the mean value theorem, we can write
\begin{eqnarray*}
&&  \Big|  \sum_{j=0}^{n-1} |u_{t_{j}^{n}}|^{\frac{1}{H}}(t_{j+1}^{n} -t_{j}^{n}) -\dint_0^T|u_s|^{\frac{1}{H}} ds\Big|
  \leq   \sum_{j=0}^{n-1}\int_{t^n_j}^{t^n_{j+1}}\left||u_{t^n_j}|^{\frac{1}{H}}-|u_s|^{\frac{1}{H}}\right| ds  \\
 &  &  \qquad  \qquad \leq   \frac 1H \sum_{j=0}^{n-1}\int_{t^n_j}^{t^n_{j+1}}|u_{t^n_j}-u_s| \left(|u_{t^n_j}|^{\frac{1}{H} -1}+|u_s|^{\frac{1}{H} -1}\right) ds.
\end{eqnarray*}
Then, applying  H\"{o}lder inequality and assumption (\ref{assump11}), yields 
\begin{eqnarray*}
\E \left( \left| \sum_{j=0}^{n-1} |u_{t_{j}^{n}}|^{\frac{1}{H}}(t_{j+1}^{n} -t_{j}^{n}) -\dint_0^T|u_s|^{\frac{1}{H}} ds\right|  \right)
  &\leq& C \sum_{j=1}^{n-1}\dint_{t^n_j}^{t^n_{j+1}} (t^n_j)^{-\alpha_1} ({t^n_{j+1}} -{t^n_{j}})^{\gamma} ds  + Cn^{-1} \\
  &\leq &C n^{\alpha_1-{\gamma}-1} \sum_{i=1}^{n-1}{i^{-{\alpha_1}}}+Cn^{-1}.
\end{eqnarray*}
This proves that $ \sum_{j=0}^{n-1} |u_{t_{j}^{n}}|^{\frac{1}{H}}(t_{j+1}^{n} -t_{j}^{n})$
converge in $L^1$ to $\dint_0^T|u_s|^{\frac{1}{H}} ds$ as $n$ tends to infinity.  This convergence, together with  (\ref{equ8}), imply that  $D_n$ converges to zero as $n$  goes to infinity, which concludes the proof of the theorem. 
\eop

\section{Divergence integral with respect to a $d$-dimensional fBm}
The purpose of this section is to generalize Theorem \ref{the2}   to multidimensional processes. In order to proceed with this generalization, we first introduce the following notation. 
 Consider a $d$-dimensional fractional Brownian motion  ($d\ge 2$)
$$
B=\{B_t, t\in [0,T]\} = \{ (B_t^{(1)}, B_t^{(2)},\dots, B_t^{(d)}),\,\, {t\in [0,T]}\}
$$
with Hurst parameter $H\in (0,1)$ defined in
a complete probability space $(\Omega, \mathcal{F},P)$, where $\mathcal{F}$ is generated by $B$. That is, the components $B^{(i)}$, $i=1,\dots,d$, are independent fractional Brownian motions with Hurst parameter $H$. We can define the derivative and divergence operators, $D^{(i)}$ and $\delta^{(i)}$, with respect to each component $B^{(i)},$ as in Section 2. Denote by $\mathbb{D}_i^{1,p}(\HH)$ the associated Sobolev spaces. We assume that these spaces include functionals depending on of all the components of $B$ and not only the $i$th component.
 
The Hilbert space $\HH_d$ associated with $B$ is the completion of the space $\mathcal{E}_d$ of step functions
$\varphi =(\varphi^{(1)},\dots,\varphi^{(d)}) : [0,T]\rightarrow \R^d$ with respect to the inner product
\[
\langle \varphi, \phi \rangle_{\HH_d}  =\sum_{k=1}^d \langle \varphi^{(k)}, \phi^{(k)} \rangle_\HH.
\]
We can develop a Malliavin calculus for the process $B$, based on the Hilbert space $\HH_d$. 
We denote by $\mathcal{S}_{d}$ the space of smooth and cylindrical random variables of the
form
 $$
 F=f\left( B(\varphi _{1}), \ldots ,
B(\varphi_{n})\right),
$$ 
where  $f\in C_{b}^{\infty}(\R^{n})$,
$\varphi_{j} =(\varphi_{j}^{(1)},\dots,\varphi_{j}^{(d)}) \in \mathcal{E}_d$, and   $B(\varphi_{j}) =\displaystyle\sum_{k=1}^{d} B ^{(k)} (\varphi_{j}^{(k)})$.

Denote by $\langle  \cdot ,  \cdot \rangle$ the usual inner product on $\R^d$. The following result  has been proved in \cite{GN} using the Ergodic Theorem.
\begin{lemma} \label{lem3}
Let $F$ be a bounded random variable with values in $\R^d$.  Then, we have 
$$
V_n^{\frac{1}{H}}(\langle F, B \rangle) \overset{L^{1}(\Omega)}{\longrightarrow}  \dint_{\R^d}\left[\dint_  0^T|\langle F, \xi\rangle|^{\frac{1}{H}}ds\right]\nu(d\xi),
$$
as $n$ tends to infinity, where $\nu$ is the normal distribution $N(0, I)$ on $\R^d$.
\end{lemma}
The following theorem is the multidimensional version of Theorem \ref{the2}. 
\begin{theorem}\label{the5}
Suppose that for each $i=1,\dots, d$, $u^{(i)}\in \mathbb{D}^{1, 2}(\HH)$ is a stochastic process  satisfying Hypothesis $\bf{(A.3)}$. Set $u_t=(u_t^{(1)},\dots,u_t^{(d)})$ and consider the divergence integral process $X=\{X_t, t \in [0,T]\}$ defined by $X_t :=\sum_{i=1}^d \int_0^t u_s^{(i)}\delta B_s^{(i)}$.  Then, we have 
$$
V_n^{\frac{1}{H}}(X) \overset{L^{1}(\Omega)}{\longrightarrow} \dint_{\R^d}\left[\dint_  0^T|\langle u_s, \xi\rangle|^{\frac{1}{H}}ds\right]\nu(d\xi),
$$
as $n$ tends to infinity, where $\nu$ is the normal distribution $N(0, I)$ on $\R^d$.
\end{theorem}
\bop. This theorem can be proved by the same arguments as in the  proof of Theorem \ref{the2}. We need to show that the expression
\[
F_n:= \E\left(\left|\sum_{i=0}^{n-1}\left|\sum_{k=1}^d\dint_{t_i^n}^{t_{i+1}^n} u_s^{(k)} \delta B_s^{(k)}\right|^{\frac{1}{H}}-\dint_{\R^d}\bigg[\dint_  0^T|\langle u_s, \xi\rangle|^{\frac{1}{H}}ds\bigg]\nu(d\xi)\right|\right),
\]
converges to zero as $n$ tends to infinity. 
Using the decomposition (\ref{decom}) for $\dint_{t_i^n}^{t_{i+1}^n} u_s^{(k)} \delta B_s^{(k)}$, and applying the same techniques as in the proof of Theorem \ref{the2}, it is not difficult to see that 
\begin{equation*}\label{Fn}
F_n  \le CA_{n}^H(B_n + C_n)^{1-H} + D_n,
\end{equation*}
where $B_n,$ $C_n$ are bounded, $A_n$ converges to zero as $n$ tends to infinity, and $D_n$ is given by 
\begin{eqnarray*}
 D_n &:=& \E\left(\left|\sum_{i=0}^{n-1}| \langle u_{t_{i}^{n}}, \Delta_{i}^{n}B\rangle |^{\frac{1}{H}}-\dint_{\R^d}\bigg[\dint_  0^T|\langle u_s, \xi\rangle|^{\frac{1}{H}}ds\bigg]\nu(d\xi)\right|\right).
\end{eqnarray*}
It only  remains to show that $D_n$ converges to zero as $n$ tends to infinity. To do this,  as in the proof of Theorem \ref{the2}, we introduce   the partition of interval $[0,T]$ given by $0=t_0^{nm}<\cdots <t_{nm}^{nm} = T$,  and we write
\begin{eqnarray} 
V^{n,m} &:=&   \notag
\left|  \sum_{i=0}^{nm-1} |\langle u_{t_{i}^{nm}}, \Delta_{i}^{nm}B\rangle|^{\frac{1}{H}}-\sum_{j=0}^{n-1}\dint_{\R^d}
 |\langle u_{t_{j}^{n}}, \xi\rangle|^{\frac{1}{H}}(t_{j+1}^{n} -t_{j}^{n})\nu(d\xi)\right|  \\  \notag
  & \leq  &  \sum_{j=0}^{n-1} \sum_{i=jm}^{(j+1)m -1} |\langle u_{t_{i}^{nm}}-u_{t_{j}^{n}}, \Delta_{i}^{nm}B\rangle|^{\frac{1}{H}}\\   \notag
  & & \qquad +\sum_{j=0}^{n-1}\bigg|\sum_{i=jm}^{(j+1)m -1}|\langle u_{t_{j}^{n}}, \Delta_{i}^{nm}B\rangle|^{\frac{1}{H}}-\dint_{\R^d}|\langle u_{t_{j}^{n}}, \xi\rangle|^{\frac{1}{H}}(t_{j+1}^{n} -t_{j}^{n})\nu(d\xi)\bigg|   \\  \notag
  & := & V_1^{n,m} + V_2^{n,m}.
\end{eqnarray}
Then, using the same arguments as in Theorem \ref{the2}, we have 
\begin{equation}\label{equ6m}
\lim_{n\rightarrow \infty}\sup_{m\ge 1}\E (V^{n,m}_1 )= 0.  
\end{equation}
On the other hand, Lemma \ref{lem3} implies that for all $n\geq 1$
\begin{equation}\label{equ7m}
\lim_{m\rightarrow \infty} \E (V^{n,m}_2 )= 0.  
\end{equation}
Moreover, it is not difficult to show that 
$$
\displaystyle\lim_{n\rightarrow \infty}\E\left|\sum_{j=0}^{n-1}\dint_{\R^d}
 |\langle u_{t_{j}^{n}}, \xi\rangle|^{\frac{1}{H}}(t_{j+1}^{n} -t_{j}^{n})\nu(d\xi)-\dint_{\R^d}\bigg[\dint_  0^T|\langle u_s, \xi\rangle|^{\frac{1}{H}}ds\bigg]\nu(d\xi)\right| =0.
 $$
 Finally, this convergence, together with  (\ref{equ6m}) and (\ref{equ7m}),  imply that $D_n$ converges to zero as $n$ tends to infinity. This completes the proof of
 Theorem \ref{the5}.\eop
\section{Fractional Bessel process}
In this section, we are going to apply the results of the previous section to the fractional Bessel process. 
Let $B$ be a $d$-dimensional fractional Brownian motion  ($d\ge 2$). 
The process  $R= \{R_t, t\in [0,T]\}$, defined by   $R_t= \|B_t\|$, 
is called the fractional Bessel process of dimension $d$ and Hurst parameter $H$.
It has been proved in \cite{CN} that, for $H> \frac12$, the fractional Bessel process $R$ has the following  representation
\begin{equation}\label{rep}
R_t =  \displaystyle\sum_{i=1}^{d}\dint_ 0^t\dfrac{B_s^{(i)}}{R_s}\delta B_s^{(i)} + H(d-1)\dint_ 0^t \dfrac{s^{2H-1}}{R_s}ds.
\end{equation}
This representation (\ref{rep}) is similar the one obtained for Bessel processes with respect to standard Brownian motion (see, for instance, Karatzas and Shreve \cite{KS}). Indeed, if $W$ is a $d$-Brownian motion and $R_t =\|W_t\|$, then 
$$
R_t =  \displaystyle\sum_{i=1}^{d}\dint_ 0^t\dfrac{W_s^{(i)}}{R_s}dW_s^{(i)} + \frac{d-1}2\dint_ 0^t \dfrac{ds}{R_s}.
$$

The goal of this section is to extend the integral representation  (\ref{rep}) to the case $H<\frac 12$. We cannot apply directly the It\^o's formula because the function $\|x\|$ is not smooth at the origin. We need the following extension of the domain of the divergence operator to processes with trajectories in  $L^{\beta}([0,T], \R^d)$, where $\beta >\frac 1{2H}$.
\begin{definition}\label{def3}
Fix $\beta >\frac 1{2H}$.
We say that a  $d$-dimensional stochastic process $u=(u^{(1)},\dots , u^{(d)})\in L^1(\Omega; L^{\beta}([0,T], \R^d))$ belongs to the extended domain of the divergence  ${\rm Dom}^*\delta$, if there exists $q>1$ such that 
\begin{equation}  \label{78}
|\E\langle u,  DF\rangle_{\HH_d}|= \left |\sum_{i=1}^{d}\E(\langle u^{(i)}, D^{(i)} F\rangle_{\HH})\right | \leq c_u \| F\|_{L^q(\Omega)}, 
\end{equation}
for every smooth  and cylindrical random variable $F \in \mathcal{S}_{d}$,  where $c_u$ is some constant depending on $u$. In this case $\delta(u)\in L^{p}(\Omega)$, where $p$ is the conjugate of $q$, is defined by the duality relationship 
$$
\E(\langle u,  DF\rangle_\HH )=\E(\delta(u) F),
$$
for every smooth and cylindrical  random variable $F \in \mathcal{S}_{d}$.
\end{definition}
Notice that the inner product  in (\ref{78}) is  well defined by formula  (\ref{ext}).  If $u\mathbf{1}_{[0,t]}$ belongs to the extended domain of the divergence, we will make use of the notation
\[
\delta(u\mathbf{1}_{[0,t]})  =\sum_{i=1}^d \int_0^t u^{(i)}_s \delta B_s^{(i)}.
\]
\begin{remark}  Notice that, since $\beta >\frac 1{2H}$, we have $\HH_{d}\subset L^{\beta}([0,T], \R^d))$ and then ${\rm Dom}\, \delta \subset {\rm Dom}^*\delta$.
\end{remark}

\begin{remark}  \label{rem5.1}
It should be noted that the process $R$ satisfies the following  
\begin{equation}\label{eq6}
\E(R_t^{-q}) =C t^{-Hq} \dint_0^{\infty} y^{d-1-q} e^{\frac{-y^2}{2}} dy:= K_q t^{-Hq},
\end{equation}
for every $q<d$, where $K_q$ is a positive constant. This property will be used later.
\end{remark}

We recall the following multidimensional It\^{o} formula for the fBm (see \cite{HMS}).
This formula requires a notion of extended domain of the divergence operator,  ${\rm Dom}^{E} \delta$  introduced in  \cite[Definition 3.9]{HMS}, which is slightly different from Definition
\ref{def3}, because we require   $u\in L^1(\Omega; L^{\beta}([0,T], \R^d))$ (instead of $u\in L^2(\Omega \times [0,T]; \mathbb{R}^d)$ and the extended divergence belongs to $L^p(\Omega)$ (instead of  $L^2(\Omega)$).  Our notion of extended domain will be useful to handle the case of the fractional Bessel process. Moreover, the class of test functionals is not the same, although this is not relevant because both classes are dense in $L^p(\Omega)$.

\begin{theorem} 
Let $B$ a $d$-dimensional fractional Brownian motion with Hurst parameter $H<\frac 12$. Suppose that $F\in C^2(\R^d)$ satisfies the growth condition
\begin{equation}\label{growth}
 \max_{x \in\R^d}\left\{|F(x)|, \left  \|\frac{\partial F}{\partial x_i}(x)\right\|, \left\| \frac{\partial^2 F}{\partial x^2_i}(x) \right\| , i=1,\dots,d\right\}\leq ce^{\lambda x^2},
\end{equation}
where $c$ and $\lambda$ are positive constants such that $\lambda <\dfrac{T^{-2H}}{4d }$. Then, for each $i=1,...,d$ and $t\in [0,T]$, the process $\bf{1}_{[0,t]}\dfrac{\partial F}{\partial x_i}(B_t) \in {\rm  Dom}^{E} \delta$, and the following formula holds
\begin{equation}  \label{ito}
F(B_t) = F(0)+ \sum_{i=1}^{d}\dint_ 0^t\dfrac{\partial F}{\partial x_i}(B_s)\delta B_s^{(i)}+
H \sum_{i=1}^{d}\dint_ 0^t\dfrac{\partial^2 F}{\partial x^2_i}(B_s)s^{2H-1}ds,
\end{equation}
where ${\rm Dom}^{E} \delta$ is the extended domain of the divergence operator  in the sense of Definition 3.9 in \cite{HMS}.
\end{theorem}

 The next result is a change of variable formula for the fractional Bessel process in the case $H<\frac{1}{2}$.
 
 \begin{theorem}\label{pro1}
Let $H<\frac{1}{2}$, and let  $R=\{R_t, \in [0,T]\}$ be the fractional Bessel process.  Set $u_t^i =\frac {B^i_t}{R_t}$  and $u_t=(u_t^{(1)},\dots,u_t^{(d)})$, for
$t\in [0,T]$. Then, we have the following results:
\begin{enumerate}
\item[(i)]  For any $t\in (0,T]$, the process $\{u_s \mathbf{1}_{[0,t]}(s), s\in [0,T]\}$ belongs to the extended domain ${\rm Dom}^*\delta$ and the representation (\ref{rep}) holds true.
\item[(ii)] If $H>\frac{1}{4}$, for any $t\in [0,T]$, the process $u\mathbf{1}_{[0,t]} $ belongs to $L^2(\Omega;\HH_d)$ and to
 the domain of $\delta$  in $L^p(\Omega)$ for any $p<d$.
\end{enumerate}
\end{theorem}
\bop.  Let us first prove part (i). Since the function $\|x\|$ is not differentiable at the origin, the It\^{o} formula (\ref{ito}) cannot be applied and we need to make a suitable  approximation. For $\varepsilon >0$, consider the function $F_{\varepsilon}(x) = (\| x\|^2 +\varepsilon^2)^{\frac12}$, which is smooth   and satisfies condition (\ref{growth}). Applying It\^{o}'s formula (\ref{ito})  we have 
\begin{equation}\label{eq7}
F_{\varepsilon}(B_t) = \varepsilon+ \sum_{i=1}^{d}\dint_ 0^t\dfrac{B_s^{(i)}}{(R_s^2 +\varepsilon^2)^{\frac12}}\delta B_s^{(i)}+
Hd \int_ 0^t \dfrac{s^{2H-1}}{(R_t^2 +\varepsilon^2)^{\frac12}}ds-H  \dint_ 0^t \frac{s^{2H-1} R_s^2}{(R_s^2 +\varepsilon^2)^{\frac32}}ds.
\end{equation}
Clearly,  $F_{\varepsilon}(B_t)$ converges to $R_t$ in $L^p$ for any $p\ge 1$.  Let $1\leq p <d$. Using Minkowski's inequality,  and taking into account Remark \ref{rem5.1}, we have 
\begin{eqnarray*}
\E\left( \left| \int_ 0^t s^{2H-1}R_s^{-1}ds\right|^p \right) 
 &\leq &\left(\int_ 0^t s^{2H-1}(\E(R_s^{-p})^{\frac{1}{p}}ds\right)^p \\
&\leq  & K_p \left(\int_ 0^t s^{-H} s^{2H-1}ds\right)^p  \le  K_p H^{-p} t^{pH}.
\end{eqnarray*}
Since for every $\varepsilon >0,$  $ \frac{s^{2H-1}}{(R_s^2 +\varepsilon^2)^{\frac12}}\leq s^{2H-1}R_s^{-1}$, the dominated convergence theorem leads to the fact that $\int_ 0^t \frac{s^{2H-1}}{(R_s^2 +\varepsilon^2)^{\frac12}}ds$  converges to $\int_ 0^t \frac{s^{2H-1}}{R_s}ds$ in $L^{p}$ for any  $1\leq p<d$, as $\varepsilon$ converges to zero. 
In the same way, we prove that $\int_ 0^t \frac{s^{2H-1}R_s^2}{(R_s^2 +\varepsilon^2)^{\frac32}}ds$  converges to $\int_ 0^t \frac{s^{2H-1}}{R_s}ds$ in $L^{p}$ for any $1\leq p<d$, as $\varepsilon$ converges to zero.
Coming back to (\ref{eq7}),  we deduce that $ \sum_{i=1}^{d}\int_ 0^t\frac{B_s^{(i)}}{(R_t^2 +\varepsilon^2)^{\frac12}}\delta B_s^{(i)}$ converges in $L^{p}$ for  any $1\leq p<d$, to some limit $G_t$, as $\varepsilon$ tends to zero. 

We are going to show that the process $u\mathbf{1}_{[0,t]}$ belongs to the extended domain of the divergence and  $\delta(u\mathbf{1}_{[0,t]})=G_t$. Let $F$ be a smooth and cylindrical random variable in $\mathcal{S}_{d}$.
For $i=1,\dots,d$,  let $u_s^{\varepsilon, (i)} =\frac{B_s^{(i)}}{(R_t^2 +\varepsilon^2)^{\frac12}}$, and $u_s^\varepsilon= (u_s^{\varepsilon, (1)}, \dots, u_s^{\varepsilon, (d)})$.     By the duality relationship we obtain
\begin{equation*}\label{duality}
\E (\langle u^{\varepsilon}\bf{1}_{[0,t]}, DF\rangle_{\HH_d}  ) =  \E(\delta(u^{\varepsilon}\bf{1}_{[0,t]}) F).
\end{equation*}
 Taking into account that  $\delta(u^{\varepsilon}\bf{1}_{[0,t]})$ converges to $G_t$ in $L^p$, and that 
\[
\lim_{\varepsilon \rightarrow 0} \E(\langle u^{\varepsilon}\bf{1}_{[0,t]}, DF\rangle_{\HH_d}) =\E(\langle u \bf{1}_{[0,t]}, DF\rangle_{\HH_d}),
\]
since the components of $u$ are bounded by one, we deduce that 
\[
\E(\langle u1_{[0,t]}, DF\rangle_{\HH_d})=\E(G_tF).
\]
This implies that  $u\mathbf{1}_{[0,t]}$ belongs to the extended domain of the divergence and $\delta(u\mathbf{1}_{[0,t]})=G_t$.

 To show part (ii), let us assume that $H>\frac14$.   We first  show that for any $i=1,\dots,d$, $u^{(i)} \in L^2( \Omega; \HH)$. We can write
\begin{eqnarray*}
|u_t^{(i)}-u_s^{(i)}| 
 &  \leq  &{|B_t^{(i)} - B_s^{(i)}|}{R_t^{-1}} + {|R_s - R_t||B_s^{(i)}|}{R_t^{-1}R_s^{-1}}   \\
&  \leq &  {\| B_t - B_s\| }{R_t^{-1}} + {|R_s - R_t|\| B_s\|}{R_t^{-1}R_s^{-1}}    \\
 & \leq & 2{\| B_t - B_s\| }{R_t^{-1}},
\end{eqnarray*}
where we have used the fact that 
$$
|R_s - R_t|=\left| \| B_t\| - \| B_s\|\right|  \leq \| B_t - B_s\|.
$$
Since $|u_t^{(i)}-u_s^{(i)}|\leq 2$, we obtain
\[
|u_t^{(i)}-u_s^{(i)}| \leq 2\left({\| B_t - B_s\| }{R_t^{-1}}\wedge 1\right),
\]
which implies
\begin{equation}\label{eqqq1}
|u_t^{(i)}-u_s^{(i)}| \leq 2\| B_t - B_s\|^{\alpha} {R_t^{-\alpha}}, 
\end{equation}
for every $\alpha\in[0,1]$.
We can write, using (\ref{est01}),
\begin{eqnarray*}
 \E(\| u_t^{(i)}\|_{\HH}^2)  & \leq & k_H  \E \left( \int_ 0^T (u_s^{(i)})^2[(T-s)^{2H-1}+ s^{2H-1}]ds  \right) \\
&& + k_H \E\left(\dint_0^T\left(\dint_s^T |u_t^{(i)}-u_s^{(i)}|(t-s)^{H-\frac32}dt\right)^2 ds\right) \\
 & :=&  k_H[N_1 + N_2].
\end{eqnarray*}
Since $|u_t^i|\leq 1$, it is clear that $N_1$ is bounded. To estimate $N_2$, choose $\alpha$, $q$ and $p$ such that $\frac{1}{2H} -1<\alpha \leq 1 $,  $1<q<\frac{d}{2\alpha}$, and $\frac{1}{p}+\frac{1}{q} =1$. Using inequality (\ref{eqqq1}) and Minkowski and H\"{o}lder inequalities, we get

\begin{eqnarray*}
N_2& \leq & 2 
\int_0^T\E\left(\dint_s^T \| B_t - B_s\|^{\alpha} {R_t^{-\alpha}}(t-s)^{H-\frac32}dt\right)^2 ds \\
 & \leq & 2 
\int_0^T\left(\int_s^T\left[ E( \| B_t -B_s\|^{2\alpha p})\right]^{\frac{1}{2p}} \left[ \E (R_t^{-2\alpha q}) \right]^{\frac{1}{2q}}(t-s)^{H-\frac32}dt\right)^2 ds \\
& \leq & 
C\int_0^T\left(\int_s^T(t-s)^{\alpha H} t^{-\alpha H}(t-s)^{ H -\frac32}dt\right)^2 ds  \\
& \leq & 
C\dint_0^T s^{-2\alpha H}(T-s)^{2(\alpha +1)H -1} ds \\
 & =& C T^{2H}\beta(-2\alpha H+1, 2(\alpha +1)H). 
\end{eqnarray*}
Hence, for $i=1,\dots,d$, $\E (\| u_t^{(i)}\|_{\HH}^2) <\infty$ and, therefore,  $u\in L^2(\Omega, \HH_d)$.  Moreover,  by the first assertion, it follows that for every $F\in \mathcal{S}_{d}$ and for $p<d$,
\begin{equation*}
 \E(\langle D F,u \mathbf{1}_{[0,t]}\rangle_{\HH_d})= \E(G_t F) \leq \|G_t\|_{p} \|F\|_{q}.
\end{equation*}
 Therefore, $u\mathbf{1}_{[0,t]}$ belongs to the domain of $\delta$ in $L^p(\Omega)$.  
\eop

Notice that, if $d>2$, then we can take $p=2$ in part (ii), and $u\mathbf{1}_{[0,t]}$ belongs to ${\rm Dom}\, \delta$.
Also, we remark that  although $u\mathbf{1}_{[0,t]}$ belongs to the (extended) domain of the divergence, this does not imply that each component $u^{(i)}\mathbf{1}_{[0,t]}$ belongs to the domain of $\delta^{(i)}$. In the next theorem, we show that under the stronger condition $2dH^2 >1$, each process $u^{(i)}$ belongs to $\mathbb{D}^{1,2}_i (\HH)$, and satisfy the Hypothesis {\bf{(A.3)} of Section 4.
 
\begin{theorem}\label{pro2}
 Suppose that  $2dH^2>1$. Let $R=\{R_t,  t\in [0,T]\}$ be the fractional Bessel process. Then, for $i=1,2,\dots,d$,  the process $u_t^{(i)}=\dfrac{B_t^{(i)}}{R_t}$ satisfies Hypothesis $\bf{(A.3)}$.
\end{theorem}
\bop. 
Fix $i=1,\dots, d$. The random variable $u_t^{(i)}$ is bounded and so, it is bounded in $L^{q}(\Omega)$   for all $q>\frac{1}{H}$. The Malliavin derivative  $D^{(i)} u^{(i)}$ is given by 
$$
D^{(i)}_su_t^{(i)} = \left(-R_t^{-3} (B_t^{(i)})^2+R_t^{-1}\right)  \bf{1}_{[0,t]}(s):= \phi_t \bf{1}_{[0,t]}(s).
$$
Notice that 
\begin{equation*}   \label{89}
\| D^{(i)}u^{(i)}_t \|_{\HH} \le 2R_t^{-1} t^{H}.
\end{equation*}
This implies  $D^{(i)}u^{(i)}_t$ is  bounded in $L^{\frac{1}{H}}(\Omega; \HH)$ because $dH>1$. Indeed, we have
 \begin{equation*}
\| D^{(i)}u^{(i)}_t\|_{L^{\frac{1}{H}}(\Omega; \HH)}  \le 2  \left(\E[R_t^{-\frac{1}{H}}]\right)^{H} t^{H} \le C.
\end{equation*}
Let us now prove that $u^{(i)}$ satisfies the inequalities (\ref{assump11}) and (\ref{assump21}), with $p=\frac{1}{H}$. 
Let $0<s\le t \le T$. 
Using estimate   (\ref{eqqq1}) and  choosing $\frac{1}{2H} -1<\alpha < Hd\wedge 1$, it follows that for $1< q< \frac{Hd}{\alpha}$ and $p_1>1$ such that $\frac{1}{p_1}+\frac{1}{q} =1$,
\begin{equation*} 
\| u_t^{(i)}-u_s^{(i)}\|_{L^{\frac{1}{H}}(\Omega)} \leq 2\left[ \E  \left( \| B_t - B_s\|^{\frac{\alpha p_1}{H}} \right) \right]^{\frac{H}{p_1}} \left(\E (R_t^{-\frac{\alpha q}{H}})  \right)^{\frac{H}{q}}
 \leq C (t-s)^{\alpha H} s^{-\alpha H}.
\end{equation*}
Hence inequality (\ref{assump11}) is satisfied with $\alpha_1 =\alpha H<\frac12$ and $\gamma =\alpha H>\frac12 -H$. 
In order to show  inequality (\ref{assump21}) with $p=\frac{1}{H}$, we first write for $0<r \le t \le T$,
\begin{eqnarray}   \notag
 \| \phi_t\bf{1}_{[0,t]}-\phi_r\bf{1}_{[0,r]} \|_{\HH}
&\leq  & \| \phi_t(\bf{1}_{[0,t]}-\bf{1}_{[0,r]}) \|_{\HH}+\| (\phi_t-\phi_r)\bf{1}_{[0,r]} \|_{\HH} \\  \notag
 & =& |\phi_t| \| \bf{1}_{[0,t]}-\bf{1}_{[0,r]} \|_{\HH}+|\phi_t-\phi_r|\| \bf{1}_{[0,r]} \|_{\HH} \\
& \leq & C\left( R_t^{-1}(t-r)^{H}+|\phi_t-\phi_r|r^{H}\right). \label{phi}
\end{eqnarray}
We have
\begin{eqnarray*}
 |\phi_t-\phi_r|& \le  &  \left|  R_t^{-3} (B^{(i)}_t)^2 -R_r^{-3} (B^{(i)}_r)^2 \right| + | R_t^{-1} - R_r^{-1} | \\
 &\le&  R_t^{-3} R_r^{-3} \left( |R_t^3-R_r^3| (B^{(i)}_r)^2  + R_t^3  |(B^{(i)}_t)^2-(B^{(i)}_r)^2 | \right)+ R_t^{-1} R_r^{-1} |R_t-R_r| \\
 &\le & \| B_t- B_r\| \left( 2R_t^{-1} R_r^{-1}+ 2R_t^{-3} R_r + R_t^{-2} + R_r^{-2} \right),
\end{eqnarray*}
and 
$$
 |\phi_t-\phi_r|\leq   |\phi_t|+  |\phi_r| \leq  2(R^{-1}_t +R^{-1}_r).
$$
Put $R_{tr} := R_t^{-1} R_r^{-1}+ R_t^{-3} R_r + R_t^{-2} + R_r^{-2}$. Then, the above inequalities imply
\begin{equation*}
 |\phi_t-\phi_r|\leq 4\left[ \left(\| B_t -B_r\|R_{tr}\right)\wedge \left(R_t^{-1}\vee R_r^{-1}\right)\right].
\end{equation*}
By using the same argument as above one can find also that 
\begin{equation*}
 |\phi_t-\phi_r|\leq 4\left[ \left(\| B_t -B_r\|R_{rt}\right)\wedge \left(R_t^{-1}\vee R_r^{-1}\right)\right].
\end{equation*}
Therefore, for every $\alpha \in [0,1]$, we can write
\begin{eqnarray}   \notag
  |\phi_t-\phi_r|  & \leq& 4\left[ \left(\| B_t -B_r\|  (R_{tr}\wedge R_{rt})\right)\wedge\left(R_t^{-1}\vee R_r^{-1}\right)   \right]\\ \notag
  & \leq &4 \| B_t -B_r\|   ^\alpha (R_{tr}^\alpha\wedge R_{rt}^\alpha)\left(R_t^{\alpha-1}\vee R_r^{\alpha-1}\right)   \\  \label{79}
 &\le  & C\| B_t -B_r\|^{\alpha}\left(R_t^{-\alpha-1}\vee R_r^{-\alpha-1}\right).
 \end{eqnarray}
  Then, substituting (\ref{79}) into  (\ref{phi}) yields
\begin{equation*}
\| \phi_t\bf{1}_{[0,t]}-\phi_r\bf{1}_{[0,r]} \|_{\HH}
\leq C\left( R_t^{-1}(t-r)^{H}+\| B_t -B_r\|^{\alpha}\left(R_t^{-\alpha-1}\vee R_r^{-\alpha-1}\right)r^{H}\right).
\end{equation*}
Choose $\alpha$, $p_1$ and $q$ such that  $\frac{1}{2H}-1<\alpha < (Hd-1)\wedge 1$,   $1<p_1<\frac{dH}{\alpha +1}$ and $\frac{1}{p_1}+\frac{1}{q}= 1$. Then, we can write
\begin{eqnarray*}
 &&  \E \left(\| \phi_t\bf{1}_{[0,t]}-\phi_r\bf{1}_{[0,r]} \|_{\HH}^{\frac{1}{H}} \right) \\
 &&\le  C \E\left[ R_t^{-\frac{1}{H}}(t-r)+  r\| B_t -B_r\|^{ \frac{\alpha}{H}}\left(R_t^{-\frac{\alpha+1}{H}}\vee R_r^{-\frac{\alpha+1}{H}}\right)\right] \\
& & \leq  C \left[C t^{-1}(t-r)+ r\left[ \E  \left( \|B_t -B_r\|^{\frac{\alpha q}{H}} \right) \right] ^{\frac{1}{q}}\left[ \E \left( \left(R_t^{-\frac{\alpha+1}{H}}\vee R_r^{-\frac{\alpha+1}{H}}\right)^{p_1} \right)\right]^{\frac{1}{p_1}}\right] \\
& & \leq C\bigg(r^{-1}(t-r)+ r^{-\alpha }(t-r)^{\alpha }\bigg) \\
 & & \leq 2C \max(r^{-1}(t-r),r^{-\alpha }(t-r)^{\alpha }),
\end{eqnarray*}
and inequality (\ref{assump21}) is satisfied with $\alpha_2 =H$ and $\gamma =\alpha H$.
Finally, for every $s \le t$, we have 
\begin{equation*}
\|D^{(i)}_su^{(i)}_t\|_{L^{\frac{1}{H}}(\Omega)} 
\leq \left(\E (R_t^{-\frac{1}{H}} ) \right)^{H}
 = C t^{-H},
\end{equation*}
and then assumption (\ref{assump3}) is satisfied with $\alpha_3 = H$.
This ends the proof of Theorem \ref{pro2}. 
\eop

\begin{remark}
If $Hd>1$, we can show, using the same arguments as in the proof of  Theorem \ref{pro2}, 
that $u^{(i)}    \in \mathbb{D}^{1,2}_i (\HH)$, for $i=1,2,\dots,d$.
\end{remark}

We now discuss  the properties of the process $\Theta= \{\Theta_t, t\in [0,T]\}$ defined by
$$
\Theta_t:= \displaystyle\sum_{i=1}^{d}\dint_ 0^t\dfrac{B_s^{(i)}}{R_t}\delta B_s^{(i)}.
$$
By Theorem \ref{pro2}, we have that for every $i=1,\dots,d$, $u_t^{(i)}=\dfrac{B_t^{(i)}}{R_t}$ satisfies Hypothesis $\bf{(A.3)}$ if $2dH^2 >1$. Therefore,  applying Theorem \ref{the5}, we have the following corollary.
\begin{corollary}\label{cor}
Suppose that $2dH^2 >1$. Then we have the following
\begin{equation*}
\begin{array}{ll}
V_n^{\frac{1}{H}}(\Theta) \overset{L^{1}(\Omega)}{{\longrightarrow}} \dint_{\R^d}\left[\dint_  0^T \left| \left\langle\dfrac{B_s}{R_s}, \xi  \right \rangle \right |^{\frac{1}{H}}ds\right]\nu(d\xi),
\end{array}
\end{equation*}
as $n$ tends to infinity, where $\nu$ is the normal distribution $N(0, I)$ on $\R^d$.
\end{corollary}
\begin{proposition}   \label{prop7}
The process $\Theta$ is $H$-self-similar.
\end{proposition}
\bop. Let $a>0$. By the representation (\ref{rep}) and the self-similarity of fBm, we have
\begin{eqnarray}\notag
\Theta_{at} &= & R_{at} -H(d-1)\dint_O^{at} \dfrac{s^{2H-1}}{R_s}ds
\\ & \overset{d}{=}& a^HR_t -H(d-1)a^H\dint_0^t \dfrac{u^{2H-1}}{R_u}du = a^H \Theta_t,\notag
\end{eqnarray}
where the symbol $\overset{d}{=}$ means that the distributions of both processes are the same. This proves that  $\Theta$ is $H$-self-similar.
\eop
\begin{remark} 
\begin{enumerate}
\item 
Corollary \ref{cor}  and Proposition \ref{prop7} imply that the process $\Theta$ and the fBm have the same $\frac 1H$-variation,  if $2dH^2>1$, and they are both $H$-self-similar. These results generalize those proved by Guerra and Nualart in \cite{GN} in the case $H<\frac 12$.
\item Let us note that although $\Theta$ and the one-dimensional fBm are both $H$-self-similar and have the same $\frac 1H$-variation,   as it is shown in \cite{HN}, it is not a fractional Brownian motion with Hurst parameter $H$.
The proof of this fact is based on the Wiener chaos expansion. Whereas, in the classical Brownian motion case it is well known, from L\'{e}vy's characterization theorem, that the process $\Theta$ is a Brownian motion.
\end{enumerate}
\end{remark}
\bf{Acknowledgements.} This work was carried out during a stay of El Hassan Essaky at Kansas University (Lawrence, KS), as a part of Fulbright program. He would like to thank KU, especially Professor David Nualart, for warm welcome and kind hospitality.

\end{document}